\documentclass[conference]{IEEEtran}

\usepackage[american]{babel}
\usepackage[T1]{fontenc}
\usepackage[utf8]{inputenc}
\usepackage{lmodern}
\usepackage[english]{isodate}

\usepackage[ruled,norelsize]{algorithm2e}

\usepackage{amsmath}
\usepackage{amssymb}
\usepackage{amsthm}
\usepackage{bm}

\usepackage{graphicx}
\usepackage{pgf, tikz}
\usepackage{pgfplots}
\pgfplotsset{compat=1.17}
\usepgfplotslibrary{statistics}
\usepackage{placeins}

\usepackage{subfig}
\usepackage{pgfplots}

\usepackage{nameref}
\usepackage[hidelinks, unicode]{hyperref} 

\usepackage[style=ieee,minbibnames=1,maxbibnames=2]{biblatex}
\usepackage[autostyle=true]{csquotes}
\usepackage{listingsutf8}
\usepackage{orcidlink}
\usepackage{siunitx}
\usepackage{subfig}
\usepackage{xcolor}

\definecolor{sign}{HTML}{b02a2d}
\definecolor{direction}{HTML}{007900}
\definecolor{regime}{HTML}{8c399e}
\definecolor{characteristic}{HTML}{1f5dc2}
\definecolor{mantissa}{HTML}{636363}
\definecolor{error}{HTML}{BD002A}
\definecolor{cellbg}{HTML}{EDEDED}

\definecolor{posit}{HTML}{b02a2d}
\definecolor{bfloat16}{HTML}{007900}
\definecolor{takum}{HTML}{1f5dc2}
\definecolor{float}{HTML}{636363}

\makeatletter
\def\lst@makecaption{%
  \def\@captype{table}%
  \@makecaption
}
\makeatother

\begin{document}

\title{%
	Evaluation of Bfloat16, Posit, and Takum Arithmetics in Sparse Linear Solvers
}
\author{%
	\IEEEauthorblockN{%
		Laslo Hunhold\,\orcidlink{0000-0001-8059-0298}}
	\IEEEauthorblockA{%
		\textit{Parallel and Distributed Systems Group}\\
		\textit{University of Cologne}\\
		Cologne, Germany\\
		\href{mailto:Laslo Hunhold <hunhold@uni-koeln.de>}
		{\texttt{hunhold@uni-koeln.de}}
	}
	\and
	\IEEEauthorblockN{%
		James Quinlan\,\orcidlink{0000-0002-2628-1651}}
	\IEEEauthorblockA{%
		\textit{Department of Computer Science}\\
		\textit{University of Southern Maine}\\
		Portland, ME, USA\\
		\href{mailto:James Quinlan <james.quinlan@maine.edu>}
		{\texttt{james.quinlan@maine.edu}}
	}
}

\maketitle

\begin{abstract}
Solving sparse linear systems lies at the core of numerous computational applications.  Consequently, understanding the numerical performance of recently proposed alternatives to the established IEEE 754 floating-point numbers, such as bfloat16 and the tapered-precision posit and takum machine number formats, is of significant interest.  This paper examines these formats in the context of widely used solvers, namely LU, QR, and GMRES, with incomplete LU preconditioning and mixed precision iterative refinement (MPIR).  This contrasts with the prevailing emphasis on designing specialized algorithms tailored to new arithmetic formats.
\par
This paper presents an extensive and unprecedented evaluation based on the SuiteSparse Matrix Collection---a dataset of real-world matrices with diverse sizes and condition numbers.  A key contribution is the faithful reproduction of SuiteSparse's UMFPACK multifrontal LU factorization and SPQR multifrontal QR factorization for machine number formats beyond single and double-precision IEEE 754.  Tapered-precision posit and takum formats show better accuracy in direct solvers and reduced iteration counts in indirect solvers. Takum arithmetic, in particular, exhibits increased stability, even at low precision.
\end{abstract}

\section{Introduction}\label{sec:intro}
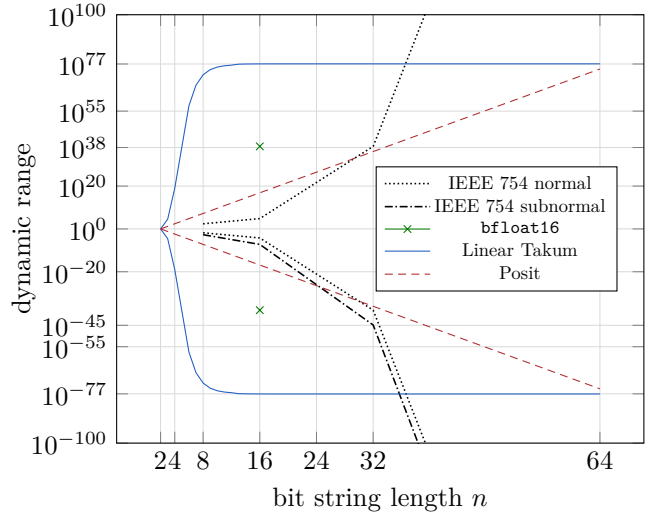
\begin{figure}[tbp]
	\begin{center}
		\begin{tikzpicture}
			\begin{axis}[
				scale only axis,
				width=\textwidth/2.6,
				height=\textwidth/3.2,
				ymin=10^-100,
				ymax=10^100,
				xlabel={bit string length $n$},
				ylabel={dynamic range},
				ymode=log,
				ylabel shift=-0.3cm,
				xtick={2,4,8,16,24,32,64},
				xminorticks=true,
				yminorticks=true,
				ytick={10^-400,10^-300,10^-200,10^-100,10^-77,10^-55,10^-45,10^-20, 1,10^20,10^38,10^55,10^77,10^100,10^200,10^300,10^400},
				grid=both,
				minor y tick num=9,
				grid style={line width=.1pt, draw=gray!10},
				major grid style={line width=.2pt,draw=gray!30},
				legend style={nodes={scale=0.7, transform shape}},
				legend style={at={(0.95,0.5)},anchor=east}
			]
				\addplot[semithick,densely dotted] table [x=n, y=ieee-normal-min, col sep=comma] {data/dynamic_range.csv};
				\addplot[forget plot,semithick,densely dotted] table [x=n, y=ieee-max, col sep=comma] {data/dynamic_range.csv};
				\addlegendentry{IEEE 754 normal};
				
				\addplot[semithick,densely dashdotted] table [x=n, y=ieee-subnormal-min, col sep=comma] {data/dynamic_range.csv};
				\addlegendentry{IEEE 754 subnormal};
				\addplot[direction,mark=x,forget plot] coordinates {(16,1.175494351e-38)};
				\addplot[direction,mark=x] coordinates {(16,3.38953139e38)};
				\addlegendentry{\texttt{bfloat16}};

				\addplot[characteristic] table [x=n, y=lintakum-min, col 
				sep=comma] {data/dynamic_range.csv};
				\addplot[characteristic,forget plot] table [x=n, y=lintakum-max, col 
				sep=comma] {data/dynamic_range.csv};
				\addlegendentry{Linear Takum};
				
				\addplot[sign,densely dashed] table [x=n, y=posit2-min, col sep=comma] {data/dynamic_range.csv};
				\addplot[sign,densely dashed,forget plot] table [x=n, y=posit2-max, col sep=comma] {data/dynamic_range.csv};
				\addlegendentry{Posit};
			\end{axis}
		\end{tikzpicture}
	\end{center}
	\caption{Dynamic range relative to the bit string length $n$ for linear takum, posit and a selection of floating-point formats.}
	\label{fig:dynamic_range}
\end{figure}
The numerical solution of sparse linear systems is a cornerstone problem in scientific computing, with applications encompassing structural analysis, circuit simulation, fluid dynamics, and machine learning. Historically, such computations have relied on the IEEE 754 floating-point standard \cite{ieee2019standard}, which has become the default format for numerical representation. However, the landscape is shifting towards low-precision formats to mitigate processor performance outpacing memory interconnect bandwidth in modern high-performance computing (commonly referred to as the \enquote{memory wall}).
\par
Emerging number formats such as bfloat16 \cite{bfloat16}, posit \cite{posits-beating_floating-point-2017}, and takum \cite{2024-takum} introduce opportunities to improve computational performance and accuracy, particularly in low-precision arithmetic. Posits and takums, for instance, employ a tapered precision scheme through variable-width exponent encoding, which allocates higher precision to values near 1 while sacrificing precision for values further from 1. Takum arithmetic represents a novel advancement over posits by offering an extensive dynamic range even at very low precisions. This design is motivated by the principle that bit-string length should primarily determine precision without imposing constraints on dynamic range---a common limitation in other formats. This paper focuses on the linear takum variant, which is a floating-point format, as opposed to the logarithmic representation in the standard takum format. For the \texttt{float8} format, which lacks standardization, we adopt the E4M3 OFP8 definition provided in \cite{ofp8} with 4 exponent bits and 3 fraction bits. Figure~\ref{fig:dynamic_range} illustrates the dynamic ranges of the formats evaluated in this study.
\par
Despite increasing interest in posits for numerical analysis \cite{posit-na-1, posit-na-2,quinlan2024iterative}, no prior work has examined takums in this context, given their novelty. Moreover, a systematic, large-scale comparison of posits and takums against bfloat16, another promising low-precision format, is currently lacking. A key question in this domain is how the expansive constant dynamic range of takums compares to that of posits, which, despite their potential, have faced adoption challenges and criticism for their limited dynamic range \cite{posits-good-bad-ugly-2019}.
\par
In this paper, we evaluate the numerical performance of these alternative number formats within a selection of established solvers that underpin scientific computing software. Our analysis is based on a large, diverse, sparse matrix test set, representing an unprecedented scale in such studies. We avoid tailoring algorithms to any specific number format to ensure unbiased results. Instead, we fully reproduce the SuiteSparse library with respect to the UMFPACK LU solver and the SPQR QR solver. This approach simulates a \enquote{blind} replacement of the underlying arithmetic in a computing environment, offering a more realistic assessment than tailored implementations. Additionally, we evaluate GMRES and mixed-precision iterative refinement methods, extending the latter to 8-bit precision---an exploration that, to the best of our knowledge, has not been previously investigated in the literature. Since most of the arithmetic formats examined in this study are implemented in software, neither computation time nor power consumption is measured; the analysis focuses solely on numerical performance.
\par
The remainder of this paper is organized as follows. Section~\ref{sec:methods} outlines the experimental methods used to benchmark the formats. Section~\ref{sec:results} presents the main results, including detailed analyses and visualizations. Finally, Section~\ref{sec:conclusion} summarizes our findings and offers conclusions.
\section{Experimental Methods}\label{sec:methods}
We evaluate the numerical performance of four fundamental approaches to solving sparse linear systems across multiple numeric formats: LU decomposition and QR factorization as direct methods, and the Generalized Minimal Residual (GMRES) method with incomplete LU preconditioning alongside Mixed Precision Iterative Refinement (MPIR) as iterative methods. These approaches encompass core techniques in numerical linear algebra, each distinguished by unique trade-offs in computational efficiency, memory consumption, and numerical stability.
\par
The benchmarking framework presented in this work, MuFoLAB (Multi-Format Linear Algebra Benchmarks) \cite{mufolab}, is designed to facilitate systematic and reproducible evaluations. It comprises three key components: a test matrix generator, a unified experimental interface for solvers, and implementations of the four solver methods under consideration. The subsequent sections provide a detailed description of these components and their roles in the benchmarking process.
\subsection{Test Matrices Generation}
\begin{figure}
	\begin{tikzpicture}
		\begin{axis}[
  				grid=major,
  				xmin=0,
  				xmax=1,
                ymax=18,
                ylabel={$\log_{10}(\kappa)$},
                ylabel shift=-0.1cm,
  				scale only axis,
   				width=\textwidth/2.6,
   				height=\textwidth/5.0,
   				xtick={0,0.25,0.5,0.75,1.0},
                xticklabel={\pgfmathparse{\tick*100}\pgfmathprintnumber{\pgfmathresult}\%},
    	       	point meta={x*100},
    		]
    		\addplot[float,thick] table [col sep=comma, x=percent, y=condition_number]{data/conditions.csv};
    	\end{axis}
    \end{tikzpicture}
    \caption{Cumulative distribution of test matrix $L^1$ condition numbers.}
    \label{fig:cond}
\end{figure}
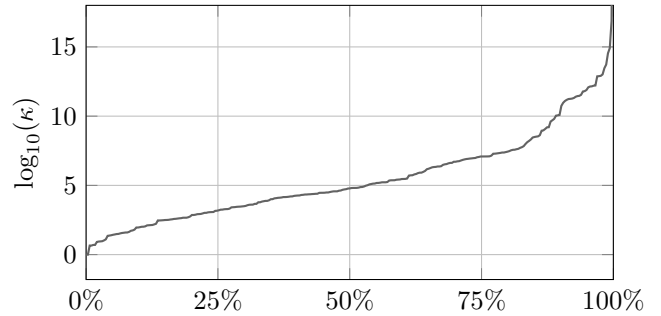

The first step in the experimental setup involves preparing a comprehensive set of sparse test matrices for benchmarking. To achieve this, we draw matrices from the SuiteSparse Matrix Collection \cite{davis2011university}, a well-established repository containing matrices from diverse application domains, including computational fluid dynamics, chemical simulation, materials science, optimal control, structural mechanics, and 2D/3D sequencing. Initially, we discard non-real matrices and those with more than $10^4$ non-zero entries, resulting in a preliminary dataset of 833 matrices. Further refinement, limited to square matrices with full rank --- criteria necessary for our benchmarks --- reduces the dataset to 295 matrices. The $L^1$ condition numbers of these matrices span several orders of magnitude, with a median value of approximately $10^5$, and around 25\% of the matrices exhibit condition numbers exceeding $10^7$, as illustrated in Figure~\ref{fig:cond}.
It should be noted that the size of the matrices is far less significant than the number of nonzero entries, as the sparse linear algebra algorithms examined in this study disregard zero entries, and any differences in fill-in affect all numerical formats equally. Furthermore, the number of nonzero entries in the matrices is limited not by inherent constraints of the number formats but to maintain reasonable computation times and reproducibility, given that the number formats are implemented in software.
\par
The established Julia package for accessing the SuiteSparse Matrix Collection (\texttt{MatrixDepot.jl}) retrieves matrices on-demand via individual internet requests. While this approach is functional in local setups with caching, it becomes unsuitable for cold-cache deployments on high-performance computing (HPC) systems and containerized environments, especially given the scale of our benchmarks, which involve hundreds of matrices. Moreover, the metadata provided by the package is frequently incomplete.
\par
To address these challenges, we introduce a second processing step to streamline and enhance matrix preparation. In this step, all filtered matrices are converted to sparse \texttt{float64} format, and complete metadata---including the number of non-zero entries, absolute minimum and maximum values, $L^1$ condition number, rank, symmetry, and positive definiteness---is computed. The processed matrices and their metadata are stored in a single compressed Julia Data File (JLD2). This self-contained file significantly improves deployment efficiency by enabling rapid access to test matrices with consistent metadata without requiring an internet connection \cite[src/TestMatricesGenerator.jl, src/TestMatrices.jl]{mufolab}.
\subsection{Common Solver Experiment Interface}
With the test matrices prepared, the next step is establishing a standardized interface for evaluating different numeric formats within a selection of solvers. Given only the matrix $A$, constructing a linear test system $Ax = b$ requires generating the solution vector $x \in \mathbb{R}^n$ and the corresponding right-hand side $b \in \mathbb{R}^n$. 
\par
Although a common approach is to set $x = (1, \dots, 1)^T$, a more representative method, as outlined in \cite[Section~5]{carson2017iterativerefinement}, involves generating $b$ randomly such that $\|b\|_\infty = 1$ and computing the reference solution by solving the system $Ax = b$ in \texttt{float128} using \texttt{Quadmath.jl} and a custom type-agnostic sparse QR solver \cite[src/QR.jl]{mufolab}. Random generation is performed using a Xoshiro pseudorandom number generator (PRNG) with a fixed seed to ensure reproducibility.
\par
For each numeric format under evaluation, the matrix $A$ and right-hand side $b$ are converted to the target type, denoted as $\tilde{A}$ and $\tilde{b}$, respectively. If any entry in $\tilde{A}$ or $\tilde{b}$ underflows or overflows during conversion, the experiment is aborted, and the failure is recorded in the results. For further discussion on dynamic range and its implications in low-precision arithmetic, refer to \cite{higham2019squeezing}. For the target types we use the Julia packages \texttt{Float8s.jl}, \texttt{BFloat16s.jl}, \texttt{Posits.jl} and \texttt{Takums.jl}.
\par
The respective solver is applied to the linear system $\tilde{A}x = \tilde{b}$, yielding an approximate solution $\tilde{x}$. This solution, cast back to \texttt{float128}, is compared to the exact solution $x$, and the absolute and relative $2$-norm errors are calculated. While it may seem counterintuitive, generating only one random sample for each matrix is sufficient because the ensemble diversity is provided by the large number of matrices in the test set. The same random seed also ensures each number format is tested with the same pseudorandom outcome \cite[src/Experiment.jl]{mufolab}.
\subsection{LU Solver}\label{sub:lu}
Although the LU decomposition algorithm is relatively straightforward to implement, the primary challenge lies in determining effective row and column permutations to minimize fill-in and enhance numerical performance. UMFPACK \cite{davis2024umfpack} is widely used and highly optimized for LU decompositions.
It employs a sophisticated rule set to select an optimal pivoting strategy, producing row and column permutations, as well as a row scaling matrix \cite[Section~1]{davis2024umfpack}.
\par
However, a significant limitation is that UMFPACK is implemented exclusively for \texttt{float32} and \texttt{float64} data types. To extend its capabilities to other numeric formats, we emulate UMFPACK's behavior by precomputing an LU decomposition in \texttt{float64} using UMFPACK for each test matrix. This computation yields the row and column permutations, which depend solely on the structural properties of the matrix and not on the specific type of its elements. The UMFPACK row scaling, however, must be computed separately for each number format, as it depends on the numerical properties of the data type. This is relatively straightforward, as it simply involves summing the absolute values in each row, yielding the vector of row 1-norms (cf. \cite[49]{davis2024umfpack}).
\par
Once the precomputed permutations and scaling factors are determined, they are applied to the matrix in the target numeric format. The system is then solved using a simple non-pivoting LU solver, effectively replicating UMFPACK's behavior. This approach ensures consistency in the decomposition process across all tested numeric formats while maintaining parity with UMFPACK's sophisticated pivoting strategy.
\subsection{QR Solver}\label{sub:qr}
Similar to the case of LU decomposition with UMFPACK, the core algorithm for QR decomposition, which employs \textsc{Householder} rotations, is relatively straightforward. However, determining optimal row and column permutations to minimize fill-in during the decomposition is the primary challenge. This challenge is addressed by SPQR, a highly optimized implementation that, like UMFPACK, is part of the SuiteSparse library \cite{davis2024spqr}. SPQR is designed to maximize numerical efficiency and reduce memory requirements through sophisticated permutation strategies.
\par
A significant limitation of SPQR is its exclusive implementation for the \texttt{float64} and \texttt{float32} data types. To enable the use of other number formats, we leverage the fact that row and column permutations are solely dependent on the structural properties of the matrix and not on the specific numeric type of its elements. Accordingly, we precompute a QR decomposition in \texttt{float64} using SPQR to extract the optimal row and column permutations. 
\par
Once the matrix is permuted according to these precomputed permutations, we collect all non-zero entries below the diagonal and apply one \textsc{Householder} rotation per column. We store both the rotation vector and the vector of indices affected by each rotation. This process effectively emulates the behavior of SPQR while allowing for a broader range of number formats, ensuring consistency and efficiency in the decomposition across all tested data types.
\subsection{Mixed Precision Iterative Refinement (MPIR) Solver}\label{sub:mpir}
The iterative refinement technique \cite{wilkinson1948progress} is a classical approach to improving an approximate solution $\tilde{x}$ to a linear system $Ax = b$. The method iteratively refines $\tilde{x}$ by solving a correction equation $Ac = r$, where $r = b - A\tilde{x}$ is the residual, and updating $\tilde{x} \leftarrow \tilde{x} + c$. This process is repeated until a convergence criterion is met, which, in our implementation, is based on the normwise backward error \cite{higham-2002-accuracy}. However, iterative refinement may fail to converge if the low-precision arithmetic causes $A$ to appear singular, thereby preventing the accurate computation of $\tilde{x}$.
\par
Mixed-precision iterative refinement (MPIR) extends this method by employing different levels of precision to optimize computational efficiency and solution accuracy \cite{carson2018accelerating}. Specifically, MPIR uses three distinct precision levels:  
\begin{enumerate}
    \item{
        \emph{Working precision} ($W$), where $A$, $b$, and
        $\tilde{x}$ are stored.
    }
    \item{
        \emph{Low precision} ($L$) for the factorization of $A$, typically to reduce computational cost.
    }
    \item{
        \emph{High precision} ($H$) for the residual calculation, ensuring accurate error correction.
    }
\end{enumerate}
These precision levels are collectively represented as a triple $(L, W, H)$.
\par
We evaluated several precision configurations, including a novel $(8, 16, 32)$ setup, alongside established configurations such as $(16, 16, 32)$, $(16, 32, 32)$, and $(16, 32, 64)$. These configurations were tested across multiple number formats, including IEEE floating-point, \texttt{bfloat16}, linear takums, and posits.
\par
For each configuration, the error tolerance was adjusted to align with the precision levels: $10^{-3}$ for $(8, 16, 32)$ and $(16, 16, 32)$, $10^{-6}$ for $(16, 32, 32)$, and $10^{-9}$ for $(16, 32, 64)$. The maximum number of iterations was set to $100$ for all experiments.
\subsection{Incomplete LU Preconditioned GMRES Solver}\label{sub:gmres}
Our implementation employs left-looking level 0 incomplete LU factorization (ILU(0)) via \texttt{IncompleteLU.jl} to reduce fill-in and use the resulting factors as a preconditioner in the Generalized Minimal Residual (GMRES) method. This combination balances computational efficiency and accelerated convergence for sparse linear systems while keeping memory requirements manageable.
\par
We leverage the \texttt{IterativeSolvers.jl} package, adhering largely to the GMRES default parameters. Specifically, we use a restart value of $\min(20, n)$, a maximum iteration count of $n$, and the modified \textsc{Gram}-\textsc{Schmidt} process for orthogonalization. For the relative tolerance, which defaults to the square root of the machine epsilon of the working precision, we instead use the square root of the machine precision of the corresponding reference float type. This approach ensures fairness across all numeric formats. For example, \texttt{float8} is used for all 8-bit types, \texttt{float16} for all 16-bit types, \texttt{float32} for all 32-bit types, and \texttt{float64} for all 64-bit types. Without this adjustment, numeric types with smaller machine epsilons than their IEEE 754 counterparts would be disproportionately disadvantaged in terms of convergence criteria.
\section{Results}\label{sec:results}
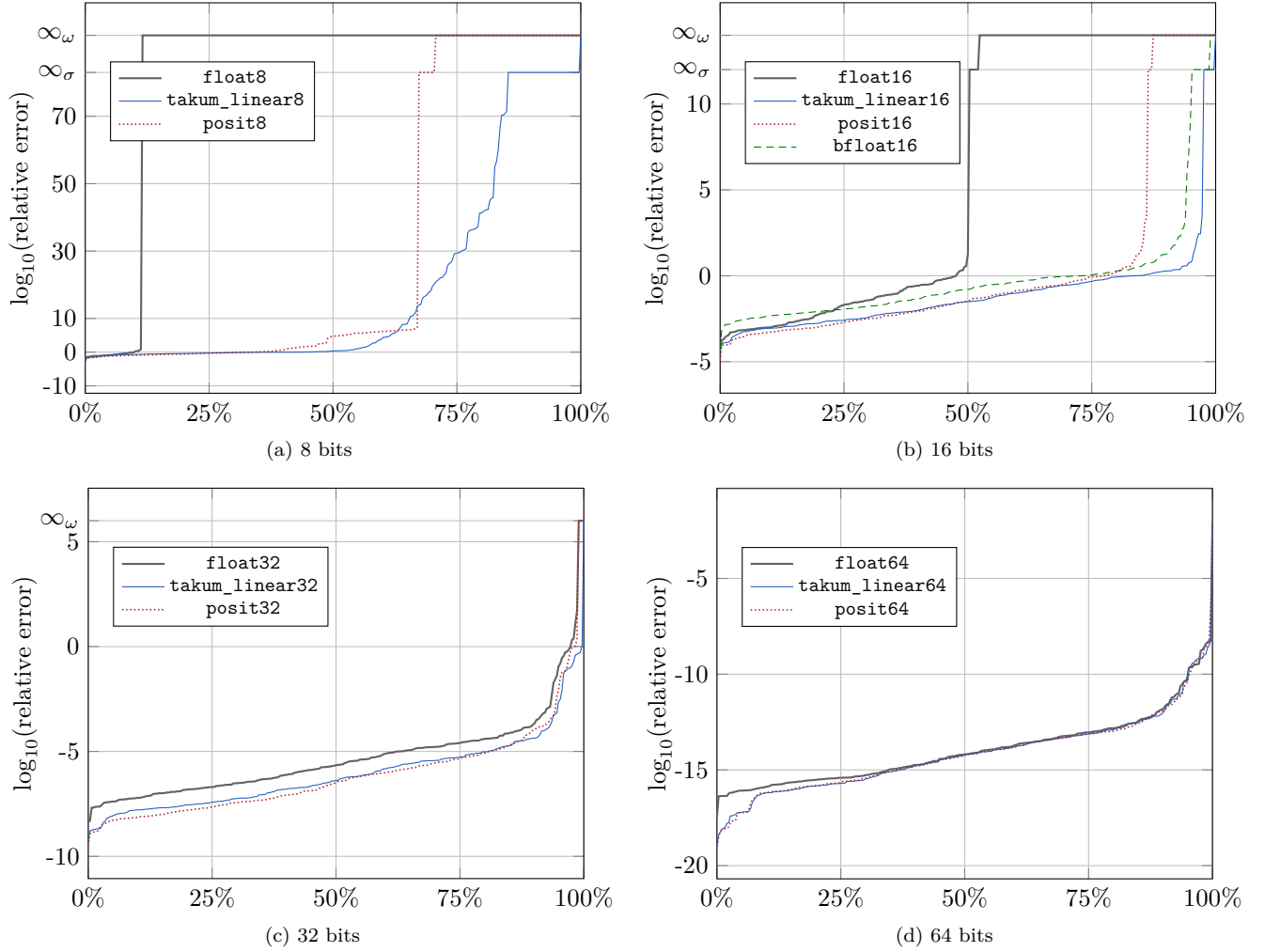
\begin{figure*}[tb] 
	\begin{center}
        \subfloat[8 bits]{
    		\begin{tikzpicture}
    			\begin{axis}[
    				scaled y ticks={base 10:-2},
    				ytick={1e-10,1e0,1e10,1e30,1e50,1e70},
    				yticklabels={-10,0,10,30,50,70},
    				extra y ticks={1e83,1e94},
    				extra y tick labels={$\infty_\sigma$, $\infty_\omega$},
    				grid=major,
    				ymode=log,
    				xmin=0,
    				xmax=1,
    				ylabel={$\log_{10}(\text{relative error})$},
                    ylabel shift=-0.2cm,
    				scale only axis,
    				width=\textwidth/2.6,
    				height=\textwidth/3.3,
    				xtick={0,0.25,0.5,0.75,1.0},
    				xticklabel={\pgfmathparse{\tick*100}\pgfmathprintnumber{\pgfmathresult}\%},
    				point meta={x*100},
    				legend style={nodes={scale=0.8, transform shape}},
    				legend style={at={(0.05,0.85)},anchor=north west},
    				]
    				\addplot[float,thick] table [col sep=comma, x=percent, y=Float8]{data/solve_lu/relative_error.sorted.csv};
    				    \addlegendentry{\texttt{float8}};
    				\addplot[takum] table [col sep=comma, x=percent, y=LinearTakum8]{data/solve_lu/relative_error.sorted.csv};
    				    \addlegendentry{\texttt{takum\_linear8}};
    				\addplot[posit,densely dotted,semithick] table [col sep=comma, x=percent, y=Posit8]{data/solve_lu/relative_error.sorted.csv};
    				    \addlegendentry{\texttt{posit8}};
    			\end{axis}
    		\end{tikzpicture}
        }
        \subfloat[16 bits]{
    		\begin{tikzpicture}
    			\begin{axis}[
    				scaled y ticks={base 10:-2},
    				ytick={1e-5,1e0,1e5,1e10},
    				yticklabels={-5,0,5,10},
    				extra y ticks={1e-7,1e12,1e14},
    				extra y tick labels={$-\infty$,$\infty_\sigma$, $\infty_\omega$},
    				grid=major,
    				ymode=log,
    				xmin=0,
    				xmax=1,
    				ylabel={$\log_{10}(\text{relative error})$},
                    ylabel shift=-0.2cm,
    				scale only axis,
    				width=\textwidth/2.6,
    				height=\textwidth/3.3,
    				xtick={0,0.25,0.5,0.75,1.0},
    				xticklabel={\pgfmathparse{\tick*100}\pgfmathprintnumber{\pgfmathresult}\%},
    				point meta={x*100},
    				legend style={nodes={scale=0.8, transform shape}},
    				legend style={at={(0.05,0.85)},anchor=north west},
    				]
    				\addplot[float,thick] table [col sep=comma, x=percent, y=Float16]{data/solve_lu/relative_error.sorted.csv};
    				    \addlegendentry{\texttt{float16}};
    				\addplot[takum] table [col sep=comma, x=percent, y=LinearTakum16]{data/solve_lu/relative_error.sorted.csv};
    				    \addlegendentry{\texttt{takum\_linear16}};
    				\addplot[posit,densely dotted,semithick] table [col sep=comma, x=percent, y=Posit16]{data/solve_lu/relative_error.sorted.csv};
    				    \addlegendentry{\texttt{posit16}};
    				\addplot[bfloat16,densely dashed] table [col sep=comma, x=percent, y=BFloat16]{data/solve_lu/relative_error.sorted.csv};
    				    \addlegendentry{\texttt{bfloat16}};
    			\end{axis}
    		\end{tikzpicture}
        }
        \\
        \subfloat[32 bits]{
    		\begin{tikzpicture}
    			\begin{axis}[
    				scaled y ticks={base 10:-2},
    				ytick={1e-10,1e-5,1e0,1e5,1e10},
    				yticklabels={-10,-5,0,5,10},
    				extra y ticks={1e6},
    				extra y tick labels={$\infty_\omega$},
    				grid=major,
    				ymode=log,
    				xmin=0,
    				xmax=1,
    				ylabel={$\log_{10}(\text{relative error})$},
                    ylabel shift=-0.2cm,
    				scale only axis,
    				width=\textwidth/2.6,
    				height=\textwidth/3.3,
    				xtick={0,0.25,0.5,0.75,1.0},
    				xticklabel={\pgfmathparse{\tick*100}\pgfmathprintnumber{\pgfmathresult}\%},
    				point meta={x*100},
    				legend style={nodes={scale=0.8, transform shape}},
    				legend style={at={(0.05,0.85)},anchor=north west},
    				]
    				\addplot[float,thick] table [col sep=comma, x=percent, y=Float32]{data/solve_lu/relative_error.sorted.csv};
    				    \addlegendentry{\texttt{float32}};
    				\addplot[takum] table [col sep=comma, x=percent, y=LinearTakum32]{data/solve_lu/relative_error.sorted.csv};
    				    \addlegendentry{\texttt{takum\_linear32}};
    				\addplot[posit,densely dotted,semithick] table [col sep=comma, x=percent, y=Posit32]{data/solve_lu/relative_error.sorted.csv};
    				    \addlegendentry{\texttt{posit32}};
    			\end{axis}
    		\end{tikzpicture}
        }
        \subfloat[64 bits]{
    		\begin{tikzpicture}
    			\begin{axis}[
    				scaled y ticks={base 10:-2},
    				ytick={1e-20,1e-15,1e-10,1e-5,1e-0},
    				yticklabels={-20,-15,-10,-5,0},
    				extra y ticks={},
    				extra y tick labels={},
    				grid=major,
    				ymode=log,
    				xmin=0,
    				xmax=1,
    				ylabel={$\log_{10}(\text{relative error})$},
                    ylabel shift=-0.2cm,
    				scale only axis,
    				width=\textwidth/2.6,
    				height=\textwidth/3.3,
    				xtick={0,0.25,0.5,0.75,1.0},
    				xticklabel={\pgfmathparse{\tick*100}\pgfmathprintnumber{\pgfmathresult}\%},
    				point meta={x*100},
    				legend style={nodes={scale=0.8, transform shape}},
    				legend style={at={(0.05,0.85)},anchor=north west},
    				]
    				\addplot[float,thick] table [col sep=comma, x=percent, y=Float64]{data/solve_lu/relative_error.sorted.csv};
    				    \addlegendentry{\texttt{float64}};
    				\addplot[takum] table [col sep=comma, x=percent, y=LinearTakum64]{data/solve_lu/relative_error.sorted.csv};
    				    \addlegendentry{\texttt{takum\_linear64}};
    				\addplot[posit,densely dotted,semithick] table [col sep=comma, x=percent, y=Posit64]{data/solve_lu/relative_error.sorted.csv};
    				    \addlegendentry{\texttt{posit64}};
    			\end{axis}
    		\end{tikzpicture}
        }
	\end{center}
	\caption{
        Cumulative error distribution of the relative errors of
        the solutions of the linear systems via fully pivoted
        LU decomposition using a range of machine number types.
        The symbol $\infty_\sigma$ denotes where the conversion
        of the matrix to the target number type turned it singular,
        $\infty_\omega$ denotes where the dynamic range of the
        matrix entries exceeded the target number type.
	}
	\label{fig:lu}
\end{figure*}

The results for the LU solver are presented in Figure~\ref{fig:lu}. As shown, both \texttt{posit8} and, to an even greater extent, \texttt{takum\_linear8} significantly outperform \texttt{float8} in terms of solution accuracy. This pattern persists across 16, 32, and 64 bits, with posits and takums consistently surpassing or at least matching the corresponding IEEE 754 floating-point types. 
\par
An especially noteworthy observation is that \texttt{takum\_linear16} consistently outperforms \texttt{bfloat16}, whereas \texttt{posit16} exhibits reduced accuracy for the lower quartile of matrices. Interestingly, while \texttt{bfloat16} generally achieves higher accuracy than \texttt{float16}, it is less accurate for approximately 25\% of the test matrices. This behavior underscores the nuanced trade-offs among these numeric formats and suggests that takums offer greater dependability in challenging cases.
\par
\begin{figure*}[tb] 
	\begin{center}
        \subfloat[8 bits]{
    		\begin{tikzpicture}
    			\begin{axis}[
    				scaled y ticks={base 10:-2},
    				ytick={1e-4,1e0,1e10,1e20,1e30},
    				yticklabels={-5,0,10,20,30},
    				extra y ticks={1e39,1e44},
    				extra y tick labels={$\infty_\sigma$, $\infty_\omega$},
    				grid=major,
    				ymode=log,
    				xmin=0,
    				xmax=1,
    				ylabel={$\log_{10}(\text{relative error})$},
                    ylabel shift=-0.2cm,
    				scale only axis,
    				width=\textwidth/2.6,
    				height=\textwidth/3.3,
    				xtick={0,0.25,0.5,0.75,1.0},
    				xticklabel={\pgfmathparse{\tick*100}\pgfmathprintnumber{\pgfmathresult}\%},
    				point meta={x*100},
    				legend style={nodes={scale=0.8, transform shape}},
    				legend style={at={(0.05,0.5)},anchor=west},
    				]
    				\addplot[float,thick] table [col sep=comma, x=percent, y=Float8]{data/solve_qr/relative_error.sorted.csv};
    				    \addlegendentry{\texttt{float8}};
    				\addplot[takum] table [col sep=comma, x=percent, y=LinearTakum8]{data/solve_qr/relative_error.sorted.csv};
    				    \addlegendentry{\texttt{takum\_linear8}};
    				\addplot[posit,densely dotted,semithick] table [col sep=comma, x=percent, y=Posit8]{data/solve_qr/relative_error.sorted.csv};
    				    \addlegendentry{\texttt{posit8}};
    			\end{axis}
    		\end{tikzpicture}
        }
        \subfloat[16 bits]{
    		\begin{tikzpicture}
    			\begin{axis}[
    				scaled y ticks={base 10:-2},
    				ytick={1e-5,1e0,1e10,1e20,1e30},
    				yticklabels={-5,0,10,20,30},
    				extra y ticks={1e42,1e48},
    				extra y tick labels={$\infty_\sigma$, $\infty_\omega$},
    				grid=major,
    				ymode=log,
    				xmin=0,
    				xmax=1,
    				ylabel={$\log_{10}(\text{relative error})$},
                    ylabel shift=-0.2cm,
    				scale only axis,
    				width=\textwidth/2.6,
    				height=\textwidth/3.3,
    				xtick={0,0.25,0.5,0.75,1.0},
    				xticklabel={\pgfmathparse{\tick*100}\pgfmathprintnumber{\pgfmathresult}\%},
    				point meta={x*100},
    				legend style={nodes={scale=0.8, transform shape}},
    				legend style={at={(0.05,0.5)},anchor=west},
    				]
    				\addplot[float,thick] table [col sep=comma, x=percent, y=Float16]{data/solve_qr/relative_error.sorted.csv};
    				    \addlegendentry{\texttt{float16}};
    				\addplot[takum] table [col sep=comma, x=percent, y=LinearTakum16]{data/solve_qr/relative_error.sorted.csv};
    				    \addlegendentry{\texttt{takum\_linear16}};
    				\addplot[posit,densely dotted,semithick] table [col sep=comma, x=percent, y=Posit16]{data/solve_qr/relative_error.sorted.csv};
    				    \addlegendentry{\texttt{posit16}};
    				\addplot[bfloat16,densely dashed] table [col sep=comma, x=percent, y=BFloat16]{data/solve_qr/relative_error.sorted.csv};
    				    \addlegendentry{\texttt{bfloat16}};
    			\end{axis}
    		\end{tikzpicture}
        }
        \\
        \subfloat[32 bits]{
    		\begin{tikzpicture}
    			\begin{axis}[
    				scaled y ticks={base 10:-2},
    				ytick={1e-10,1e-5,1e0,1e5,1e10},
    				yticklabels={-10,-5,0,5,10},
    				extra y ticks={1e6},
    				extra y tick labels={$\infty_\omega$},
    				grid=major,
    				ymode=log,
    				xmin=0,
    				xmax=1,
    				ylabel={$\log_{10}(\text{relative error})$},
                    ylabel shift=-0.2cm,
    				scale only axis,
    				width=\textwidth/2.6,
    				height=\textwidth/3.3,
    				xtick={0,0.25,0.5,0.75,1.0},
    				xticklabel={\pgfmathparse{\tick*100}\pgfmathprintnumber{\pgfmathresult}\%},
    				point meta={x*100},
    				legend style={nodes={scale=0.8, transform shape}},
    				legend style={at={(0.05,0.95)},anchor=north west},
    				]
    				\addplot[float,thick] table [col sep=comma, x=percent, y=Float32]{data/solve_qr/relative_error.sorted.csv};
    				    \addlegendentry{\texttt{float32}};
    				\addplot[takum] table [col sep=comma, x=percent, y=LinearTakum32]{data/solve_qr/relative_error.sorted.csv};
    				    \addlegendentry{\texttt{takum\_linear32}};
    				\addplot[posit,densely dotted,semithick] table [col sep=comma, x=percent, y=Posit32]{data/solve_qr/relative_error.sorted.csv};
    				    \addlegendentry{\texttt{posit32}};
    			\end{axis}
    		\end{tikzpicture}
        }
        \subfloat[64 bits]{
    		\begin{tikzpicture}
    			\begin{axis}[
    				scaled y ticks={base 10:-2},
    				ytick={1e-20,1e-10,1e0},
    				yticklabels={-20,-10,0},
    				extra y ticks={},
    				extra y tick labels={},
    				grid=major,
    				ymode=log,
    				xmin=0,
    				xmax=1,
    				ylabel={$\log_{10}(\text{relative error})$},
                    ylabel shift=-0.2cm,
    				scale only axis,
    				width=\textwidth/2.6,
    				height=\textwidth/3.3,
    				xtick={0,0.25,0.5,0.75,1.0},
    				xticklabel={\pgfmathparse{\tick*100}\pgfmathprintnumber{\pgfmathresult}\%},
    				point meta={x*100},
    				legend style={nodes={scale=0.8, transform shape}},
    				legend style={at={(0.05,0.95)},anchor=north west},
    				]
    				\addplot[float,thick] table [col sep=comma, x=percent, y=Float64]{data/solve_qr/relative_error.sorted.csv};
    				    \addlegendentry{\texttt{float64}};
    				\addplot[takum] table [col sep=comma, x=percent, y=LinearTakum64]{data/solve_qr/relative_error.sorted.csv};
    				    \addlegendentry{\texttt{takum\_linear64}};
    				\addplot[posit,densely dotted,semithick] table [col sep=comma, x=percent, y=Posit64]{data/solve_qr/relative_error.sorted.csv};
    				    \addlegendentry{\texttt{posit64}};
    			\end{axis}
    		\end{tikzpicture}
        }
	\end{center}
	\caption{
        Cumulative error distribution of the relative errors of
        the solutions of the linear systems via
        QR decomposition using a range of machine number types.
        The symbol $\infty_\sigma$ denotes where the conversion
        of the matrix to the target number type turned it singular,
        $\infty_\omega$ denotes where the dynamic range of the
        matrix entries exceeded the target number type.
	}
	\label{fig:qr}
\end{figure*}
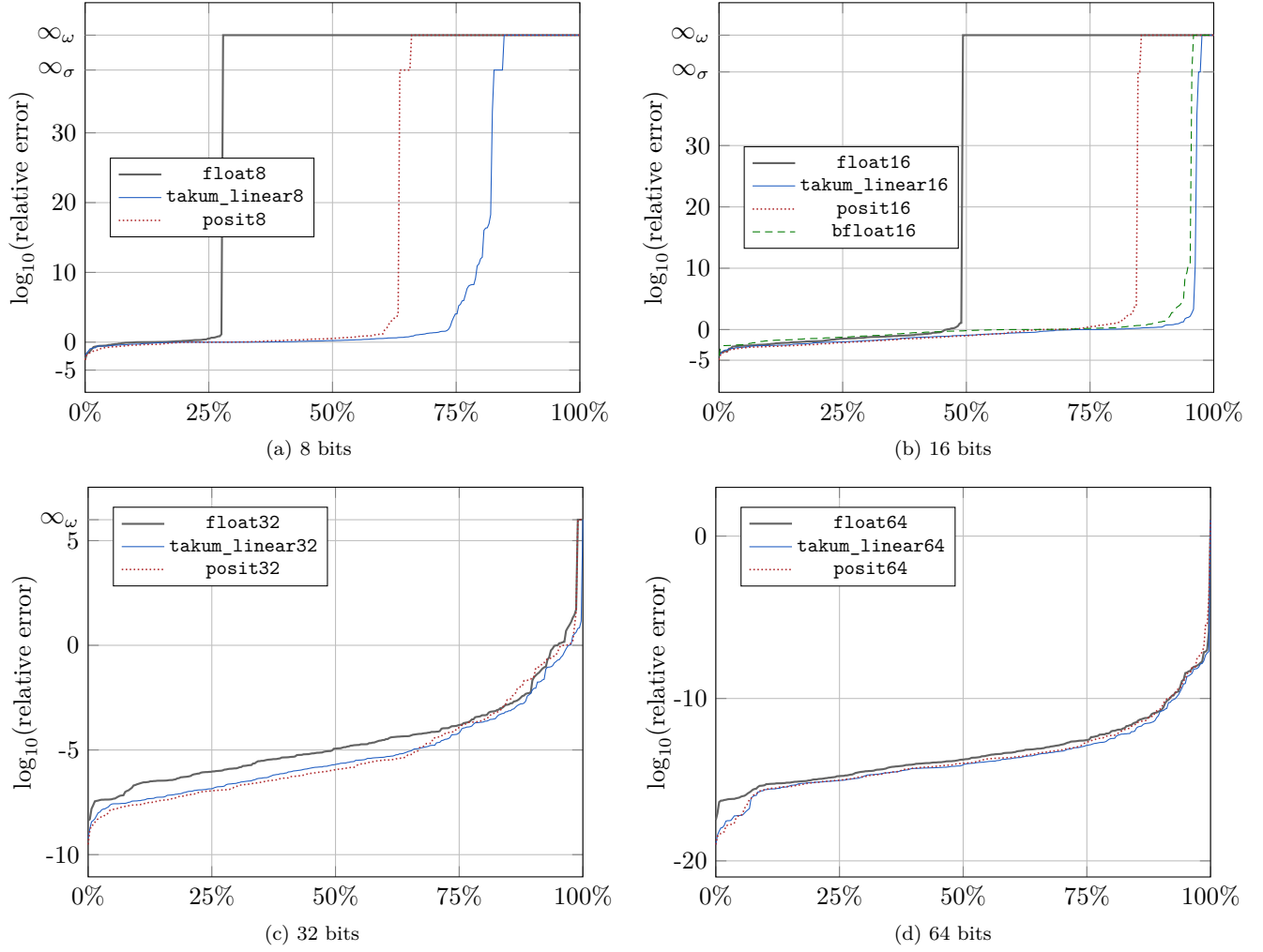

The results for the QR solver, shown in Figure~\ref{fig:qr}, exhibit a similar overall trend. Both posits and takums consistently outperform or match their corresponding IEEE 754 floating-point counterparts. Notably, \texttt{takum\_linear16} consistently surpasses \texttt{bfloat16}, showing increased accuracy across all test cases.
\par
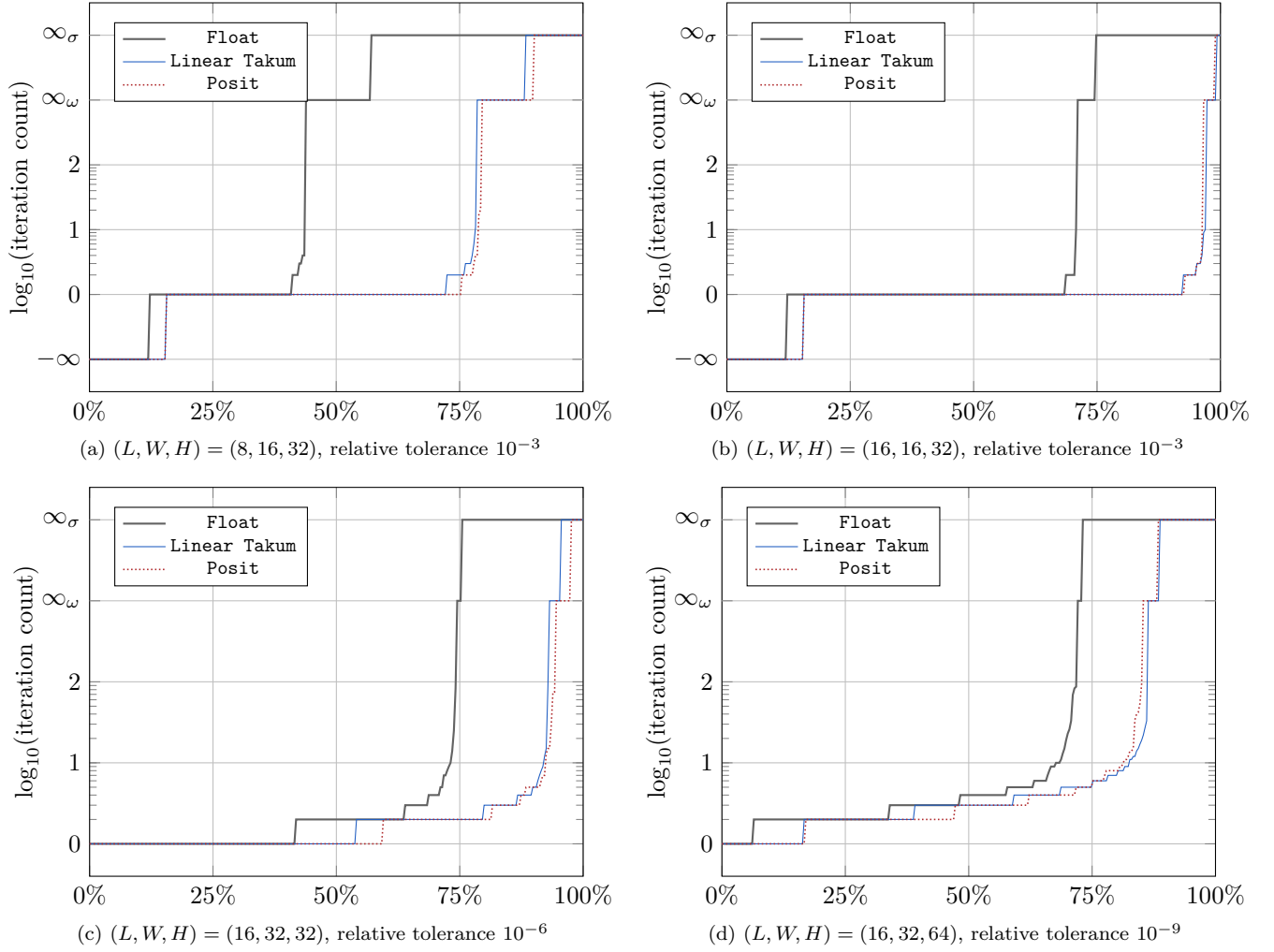
\begin{figure*}[tb] 
	\begin{center}
        \subfloat[$(L,W,H) = (8,16,32)$, relative tolerance $10^{-3}$]{
    		\begin{tikzpicture}
    			\begin{axis}[
    				scaled y ticks={base 10:-2},
    				ytick={1e0,1e1,1e2},
    				yticklabels={0,1,2},
    				extra y ticks={1e-1,1e3,1e4},
    				extra y tick labels={$-\infty$,$\infty_\omega$,$\infty_\sigma$},
    				grid=major,
    				ymode=log,
    				xmin=0,
    				xmax=1,
    				ylabel={$\log_{10}(\text{iteration count})$},
                    ylabel shift=-0.2cm,
    				scale only axis,
    				width=\textwidth/2.6,
    				height=\textwidth/3.3,
    				xtick={0,0.25,0.5,0.75,1.0},
    				xticklabel={\pgfmathparse{\tick*100}\pgfmathprintnumber{\pgfmathresult}\%},
    				point meta={x*100},
    				legend style={nodes={scale=0.8, transform shape}},
    				legend style={at={(0.05,0.95)},anchor=north west},
    				]
    				\addplot[float,thick] table [col sep=comma, x=percent, y=Float64]{data/solve_mpir_float_08_16_32/iteration_count.sorted.csv};
    				    \addlegendentry{\texttt{Float}};
    				\addplot[takum] table [col sep=comma, x=percent, y=Float64]{data/solve_mpir_takum_08_16_32/iteration_count.sorted.csv};
    				    \addlegendentry{\texttt{Linear Takum}};
    				\addplot[posit,densely dotted,semithick] table [col sep=comma, x=percent, y=Float64]{data/solve_mpir_posit_08_16_32/iteration_count.sorted.csv};
    				    \addlegendentry{\texttt{Posit}};
    			\end{axis}
    		\end{tikzpicture}
        }
        \subfloat[$(L,W,H) = (16,16,32)$, relative tolerance $10^{-3}$]{
    		\begin{tikzpicture}
    			\begin{axis}[
    				scaled y ticks={base 10:-2},
    				ytick={1e0,1e1,1e2},
    				yticklabels={0,1,2},
    				extra y ticks={1e-1,1e3,1e4},
    				extra y tick labels={$-\infty$,$\infty_\omega$,$\infty_\sigma$},
    				grid=major,
    				ymode=log,
    				xmin=0,
    				xmax=1,
    				ylabel={$\log_{10}(\text{iteration count})$},
                    ylabel shift=-0.2cm,
    				scale only axis,
    				width=\textwidth/2.6,
    				height=\textwidth/3.3,
    				xtick={0,0.25,0.5,0.75,1.0},
    				xticklabel={\pgfmathparse{\tick*100}\pgfmathprintnumber{\pgfmathresult}\%},
    				point meta={x*100},
    				legend style={nodes={scale=0.8, transform shape}},
    				legend style={at={(0.05,0.95)},anchor=north west},
    				]
    				\addplot[float,thick] table [col sep=comma, x=percent, y=Float64]{data/solve_mpir_float_16_16_32/iteration_count.sorted.csv};
    				    \addlegendentry{\texttt{Float}};
    				\addplot[takum] table [col sep=comma, x=percent, y=Float64]{data/solve_mpir_takum_16_16_32/iteration_count.sorted.csv};
    				    \addlegendentry{\texttt{Linear Takum}};
    				\addplot[posit,densely dotted,semithick] table [col sep=comma, x=percent, y=Float64]{data/solve_mpir_posit_16_16_32/iteration_count.sorted.csv};
    				    \addlegendentry{\texttt{Posit}};
    			\end{axis}
    		\end{tikzpicture}
        }
        \\
        \subfloat[$(L,W,H) = (16,32,32)$, relative tolerance $10^{-6}$]{
    		\begin{tikzpicture}
    			\begin{axis}[
    				scaled y ticks={base 10:-2},
    				ytick={1e0,1e1,1e2},
    				yticklabels={0,1,2},
    				extra y ticks={1e-1,1e3,1e4},
    				extra y tick labels={$-\infty$,$\infty_\omega$,$\infty_\sigma$},
    				grid=major,
    				ymode=log,
    				xmin=0,
    				xmax=1,
    				ylabel={$\log_{10}(\text{iteration count})$},
                    ylabel shift=-0.2cm,
    				scale only axis,
    				width=\textwidth/2.6,
    				height=\textwidth/3.3,
    				xtick={0,0.25,0.5,0.75,1.0},
    				xticklabel={\pgfmathparse{\tick*100}\pgfmathprintnumber{\pgfmathresult}\%},
    				point meta={x*100},
    				legend style={nodes={scale=0.8, transform shape}},
    				legend style={at={(0.05,0.95)},anchor=north west},
    				]
    				\addplot[float,thick] table [col sep=comma, x=percent, y=Float64]{data/solve_mpir_float_16_32_32/iteration_count.sorted.csv};
    				    \addlegendentry{\texttt{Float}};
    				\addplot[takum] table [col sep=comma, x=percent, y=Float64]{data/solve_mpir_takum_16_32_32/iteration_count.sorted.csv};
    				    \addlegendentry{\texttt{Linear Takum}};
    				\addplot[posit,densely dotted,semithick] table [col sep=comma, x=percent, y=Float64]{data/solve_mpir_posit_16_32_32/iteration_count.sorted.csv};
    				    \addlegendentry{\texttt{Posit}};
    			\end{axis}
    		\end{tikzpicture}
        }
        \subfloat[$(L,W,H) = (16,32,64)$, relative tolerance $10^{-9}$]{
    		\begin{tikzpicture}
    			\begin{axis}[
    				scaled y ticks={base 10:-2},
    				ytick={1e0,1e1,1e2},
    				yticklabels={0,1,2},
    				extra y ticks={1e-1,1e3,1e4},
    				extra y tick labels={$-\infty$,$\infty_\omega$,$\infty_\sigma$},
    				grid=major,
    				ymode=log,
    				xmin=0,
    				xmax=1,
    				ylabel={$\log_{10}(\text{iteration count})$},
                    ylabel shift=-0.2cm,
    				scale only axis,
    				width=\textwidth/2.6,
    				height=\textwidth/3.3,
    				xtick={0,0.25,0.5,0.75,1.0},
    				xticklabel={\pgfmathparse{\tick*100}\pgfmathprintnumber{\pgfmathresult}\%},
    				point meta={x*100},
    				legend style={nodes={scale=0.8, transform shape}},
    				legend style={at={(0.05,0.95)},anchor=north west},
    				]
    				\addplot[float,thick] table [col sep=comma, x=percent, y=Float64]{data/solve_mpir_float_16_32_64/iteration_count.sorted.csv};
    				    \addlegendentry{\texttt{Float}};
    				\addplot[takum] table [col sep=comma, x=percent, y=Float64]{data/solve_mpir_takum_16_32_64/iteration_count.sorted.csv};
    				    \addlegendentry{\texttt{Linear Takum}};
    				\addplot[posit,densely dotted,semithick] table [col sep=comma, x=percent, y=Float64]{data/solve_mpir_posit_16_32_64/iteration_count.sorted.csv};
    				    \addlegendentry{\texttt{Posit}};
    			\end{axis}
    		\end{tikzpicture}
        }
	\end{center}
	\caption{
        Cumulative distribution of the MPIR iteration counts
        using a range of machine number types.
        The symbol $\infty_\sigma$ denotes where the initial low-precision LU
        decomposition yielded a singular system,
        $\infty_\omega$ denotes where the maximum iteration count was
        reached without the residual going below the desired relative tolerance.
	}
	\label{fig:mpir}
\end{figure*}

The mixed precision iterative refinement results are presented in Figure~\ref{fig:mpir}. Overall, both posits and takums demonstrate significantly lower iteration counts and fewer occurrences of singularities or maximum iteration limit exceedances compared to their respective IEEE 754 floating-point counterparts. When comparing posits and takums, no clear advantage of one format emerges, as their numerical performance appears comparable across the evaluated test cases.
\par
\begin{figure*}[tb] 
	\begin{center}
        \subfloat[8 bits]{
    		\begin{tikzpicture}
    			\begin{axis}[
    				scaled y ticks={base 10:-2},
    				ytick={1e0,1e2,1e4},
    				yticklabels={0,2,4},
    				extra y ticks={1e-1,1e5},
    				extra y tick labels={$-\infty$,$\infty_\omega$},
    				grid=major,
    				ymode=log,
    				xmin=0,
    				xmax=1,
    				ylabel={$\log_{10}(\text{iteration count})$},
                    ylabel shift=-0.2cm,
    				scale only axis,
    				width=\textwidth/2.6,
    				height=\textwidth/3.3,
    				xtick={0,0.25,0.5,0.75,1.0},
    				xticklabel={\pgfmathparse{\tick*100}\pgfmathprintnumber{\pgfmathresult}\%},
    				point meta={x*100},
    				legend style={nodes={scale=0.8, transform shape}},
    				legend style={at={(0.13,0.5)},anchor=north west},
    				]
    				\addplot[float,thick] table [col sep=comma, x=percent, y=Float8]{data/solve_gmres_ilu/iteration_count.sorted.csv};
    				    \addlegendentry{\texttt{float8}};
    				\addplot[takum] table [col sep=comma, x=percent, y=LinearTakum8]{data/solve_gmres_ilu/iteration_count.sorted.csv};
    				    \addlegendentry{\texttt{takum\_linear8}};
    				\addplot[posit,densely dotted,semithick] table [col sep=comma, x=percent, y=Posit8]{data/solve_gmres_ilu/iteration_count.sorted.csv};
    				    \addlegendentry{\texttt{posit8}};
    			\end{axis}
    		\end{tikzpicture}
        }
        \subfloat[16 bits]{
    		\begin{tikzpicture}
    			\begin{axis}[
    				scaled y ticks={base 10:-2},
    				ytick={1e0,1e2,1e4},
    				yticklabels={0,2,4},
    				extra y ticks={1e-1,1e5},
    				extra y tick labels={$-\infty$,$\infty_\omega$},
    				grid=major,
    				ymode=log,
    				xmin=0,
    				xmax=1,
    				ylabel={$\log_{10}(\text{iteration count})$},
                    ylabel shift=-0.2cm,
    				scale only axis,
    				width=\textwidth/2.6,
    				height=\textwidth/3.3,
    				xtick={0,0.25,0.5,0.75,1.0},
    				xticklabel={\pgfmathparse{\tick*100}\pgfmathprintnumber{\pgfmathresult}\%},
    				point meta={x*100},
    				legend style={nodes={scale=0.8, transform shape}},
    				legend style={at={(0.05,0.95)},anchor=north west},
    				]
    				\addplot[float,thick] table [col sep=comma, x=percent, y=Float16]{data/solve_gmres_ilu/iteration_count.sorted.csv};
    				    \addlegendentry{\texttt{float16}};
    				\addplot[takum] table [col sep=comma, x=percent, y=LinearTakum16]{data/solve_gmres_ilu/iteration_count.sorted.csv};
    				    \addlegendentry{\texttt{takum\_linear16}};
    				\addplot[posit,densely dotted,semithick] table [col sep=comma, x=percent, y=Posit16]{data/solve_gmres_ilu/iteration_count.sorted.csv};
    				    \addlegendentry{\texttt{posit16}};
    				\addplot[bfloat16,densely dashed] table [col sep=comma, x=percent, y=BFloat16]{data/solve_gmres_ilu/iteration_count.sorted.csv};
    				    \addlegendentry{\texttt{bfloat16}};
    			\end{axis}
    		\end{tikzpicture}
        }
        \\
        \subfloat[32 bits]{
    		\begin{tikzpicture}
    			\begin{axis}[
    				scaled y ticks={base 10:-2},
    				ytick={1e0,1e2,1e4},
    				yticklabels={0,2,4},
    				extra y ticks={1e-1,1e5},
    				extra y tick labels={$-\infty$,$\infty_\omega$},
    				grid=major,
    				ymode=log,
    				xmin=0,
    				xmax=1,
    				ylabel={$\log_{10}(\text{iteration count})$},
                    ylabel shift=-0.2cm,
    				scale only axis,
    				width=\textwidth/2.6,
    				height=\textwidth/3.3,
    				xtick={0,0.25,0.5,0.75,1.0},
    				xticklabel={\pgfmathparse{\tick*100}\pgfmathprintnumber{\pgfmathresult}\%},
    				point meta={x*100},
    				legend style={nodes={scale=0.8, transform shape}},
    				legend style={at={(0.05,0.95)},anchor=north west},
    				]
    				\addplot[float,thick] table [col sep=comma, x=percent, y=Float32]{data/solve_gmres_ilu/iteration_count.sorted.csv};
    				    \addlegendentry{\texttt{float32}};
    				\addplot[takum] table [col sep=comma, x=percent, y=LinearTakum32]{data/solve_gmres_ilu/iteration_count.sorted.csv};
    				    \addlegendentry{\texttt{takum\_linear32}};
    				\addplot[posit,densely dotted,semithick] table [col sep=comma, x=percent, y=Posit32]{data/solve_gmres_ilu/iteration_count.sorted.csv};
    				    \addlegendentry{\texttt{posit32}};
    			\end{axis}
    		\end{tikzpicture}
        }
        \subfloat[64 bits]{
    		\begin{tikzpicture}
    			\begin{axis}[
    				scaled y ticks={base 10:-2},
    				ytick={1e0,1e2,1e4},
    				yticklabels={0,2,4},
    				extra y ticks={1e-1,1e5},
    				extra y tick labels={$-\infty$,$\infty_\omega$},
    				grid=major,
    				ymode=log,
    				xmin=0,
    				xmax=1,
    				ylabel={$\log_{10}(\text{iteration count})$},
                    ylabel shift=-0.2cm,
    				scale only axis,
    				width=\textwidth/2.6,
    				height=\textwidth/3.3,
    				xtick={0,0.25,0.5,0.75,1.0},
    				xticklabel={\pgfmathparse{\tick*100}\pgfmathprintnumber{\pgfmathresult}\%},
    				point meta={x*100},
    				legend style={nodes={scale=0.8, transform shape}},
    				legend style={at={(0.05,0.95)},anchor=north west},
    				]
    				\addplot[float,thick] table [col sep=comma, x=percent, y=Float64]{data/solve_gmres_ilu/iteration_count.sorted.csv};
    				    \addlegendentry{\texttt{float64}};
    				\addplot[takum] table [col sep=comma, x=percent, y=LinearTakum64]{data/solve_gmres_ilu/iteration_count.sorted.csv};
    				    \addlegendentry{\texttt{takum\_linear64}};
    				\addplot[posit,densely dotted,semithick] table [col sep=comma, x=percent, y=Posit64]{data/solve_gmres_ilu/iteration_count.sorted.csv};
    				    \addlegendentry{\texttt{posit64}};
    			\end{axis}
    		\end{tikzpicture}
        }
	\end{center}
	\caption{
        Cumulative distribution of the GMRES iteration counts
        using a range of machine number types.
        The symbol $\infty_\omega$ denotes where the maximum iteration
        count was exceeded without reaching the desired tolerance.
	}
	\label{fig:gmres}
\end{figure*}
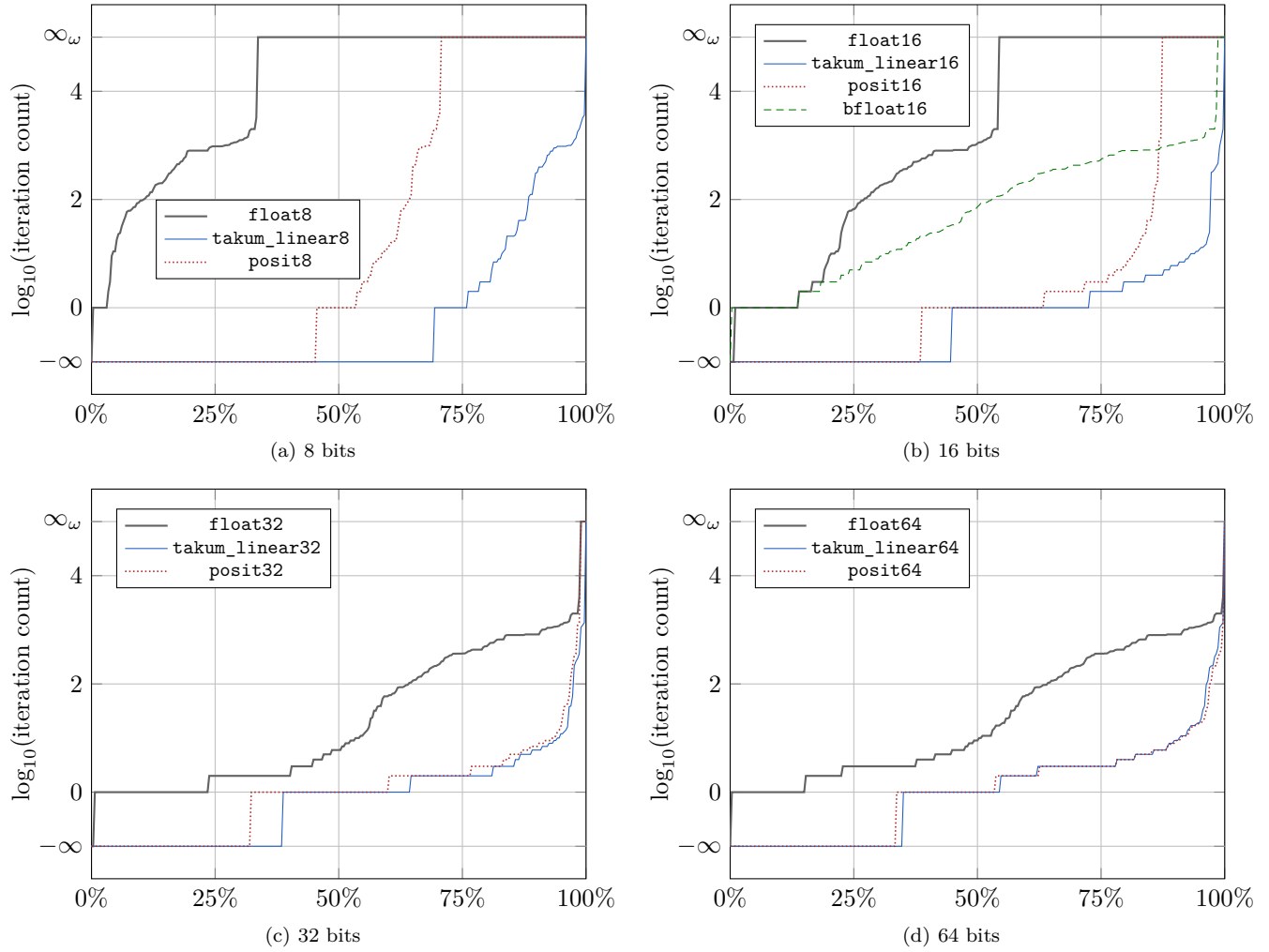

For GMRES preconditioned with incomplete LU, the results displayed in Figure~\ref{fig:gmres} highlight significant differences in performance across numeric formats. While \texttt{float8} frequently experiences overflows or requires a high number of iterations, both \texttt{posit8} and \texttt{takum\_linear8} exhibit much greater numerical stability. Notably, \texttt{takum\_linear8} avoids overflow entirely for all test matrices.
\par
This trend persists at 16 bits, where \texttt{takum\_linear16} consistently achieves lower iteration counts than \texttt{bfloat16}, in contrast to \texttt{posit16}, which occasionally lags behind. At 32 and 64 bits, posits and takums demonstrate very similar performance, both achieving significantly fewer iterations compared to \texttt{float32} and \texttt{float64}. These results underscore the advantages of tapered-precision formats in reducing computational effort and enhancing stability.
\section{Conclusion}\label{sec:conclusion}
We evaluated IEEE 754 floating-point numbers, \texttt{bfloat16}, posits, and takums across four widely used direct and iterative solving algorithms. Our experiments demonstrate that tapered-precision arithmetic consistently outperforms IEEE 754 floating-point numbers in all tested scenarios. Among the tapered-precision formats, takums exhibited exceptional performance, outperforming \texttt{bfloat16} in every case. While occasionally marginally less accurate than posits, takums delivered comparable results overall and demonstrated superior numerical stability. Notably, we successfully introduced the application of 8-bit posits and takums in mixed-precision iterative refinement, marking a possibly significant milestone in numerical computing. Additionally, GMRES exhibited particular benefit from tapered-precision formats, with takums delivering outstanding results, surpassing posits in all cases.  
\par
These findings are especially relevant as they position takums as a strong candidate to replace \texttt{bfloat16} as the state-of-the-art in 16-bit arithmetic, which posits were unable to given their limited dynamic range. Furthermore, the results address a critical question: despite having a much larger dynamic range than other number formats, including posits (see Figure~\ref{fig:dynamic_range}), takums exhibit comparable and often more favorable results. This property is possibly transformative for mixed-precision workflows, as the choice of the precision level $n$ with takums becomes purely about precision, decoupled from concerns about dynamic range.  
\par
Future research can explore further optimizations for MPIR using equilibrated matrices and investigate GMRES-based iterative refinement employing more than three precision levels \cite{amestoy2024five}.
\section*{Author Contributions}
\textbf{Laslo Hunhold}: Conceptualization, Data curation, Formal Analysis, Funding acquisition, Investigation, Methodology, Project administration, Resources, Software, Supervision, Validation, Visualization, Writing – original draft, Writing – review \& editing; \textbf{James Quinlan}: Conceptualization, Formal Analysis, Funding acquisition, Investigation, Methodology, Resources, Software, Writing – original draft, Writing – review \& editing
\printbibliography
\end{document}